\documentclass[preprint,12pt,3p]{elsarticle}

\usepackage{amssymb}
\usepackage{amsthm}

\usepackage{amsmath}
\usepackage{mathtools}
\usepackage{amsfonts}
\usepackage{mathrsfs}
\usepackage{verbatim}
\usepackage{etoolbox}
\usepackage[none]{hyphenat}

\newcommand{\R}{\mathbb{R}}

\journal{arXiv}

\begin{document}

\begin{frontmatter}

\title{An identity for the odd double factorial\tnoteref{label1}}
\tnotetext[label1]{This research did not receive any specific grant from funding agencies in the public, commercial, or not-for-profit sectors.\\ 2010 Mathematics Subject Classification 05A10}

\author{Saud Hussein}
\address{Institute of Mathematics, Academia Sinica, 6F, Astronomy-Mathematics Building, No.1, Sec.4, Roosevelt Road, Taipei 10617, Taiwan}

\ead{saudhussein@gate.sinica.edu.tw}

\begin{abstract}
The ordinary factorial may be written in terms of the Stirling numbers of the second kind as shown by Quaintance and Gould and the odd double factorial in terms of the Stirling numbers of the first kind as shown by Callan. During the preparation of an expository paper on Wolstenholme's Theorem, we discovered an expression for the odd double factorial using the Stirling numbers of the second kind. This appears to be the first such identity involving both positive and negative integers and the result is outlined here.
\end{abstract}


\end{frontmatter}


\section{Odd double factorial}

The \textit{Stirling numbers of the first kind and second kind}, denoted $s(n,k)$ and $S(n,k)$ respectively, are characterized by \begin{align}
\sum_{k=0}^n s(n,k)x^k = n!\binom{x}{n}, \quad x \in \R, \label{first} \end{align} \begin{align}\sum_{k=0}^n k!\binom{x}{k}S(n,k) = x^n, \quad x \in \R, \label{second} \end{align} and we have an explicit formula \cite[equation 13.32]{Gould} relating the two sets of numbers, \[s(n,n-k) = \sum_{j=0}^k(-1)^j\binom{n+j-1}{k+j}\binom{n+k}{k-j}S(j+k,j).\] Also, \begin{align}(-1)^ks(n,n-k) &= \sum_{j=0}^k(-1)^{j+k}\binom{n+j-1}{k+j}\binom{n+k}{k-j}S(j+k,j) \label{third}\\
&= \frac{(n+k)!}{(2k)!(n-k-1)!}\sum_{j=0}^k(-1)^{j+k}\binom{2k}{k-j}\frac{S(j+k,j)}{n+j} \nonumber\\
&= \frac{1}{(2k)!}\sum_{j=0}^k(-1)^{j+k}\binom{2k}{k-j}\frac{C(n,k)}{n+j}S(j+k,j) \nonumber
\end{align}
with $C(n,k) = (n-k)(n-k+1)\cdots (n+k)$. For fixed integer $k\geq 0$, since $(-1)^ks(n,n-k)$ is a polynomial in $n$ over the rationals with leading coefficient $\frac{(2k-1)!!}{(2k)!}$ by \cite[proposition 1.1]{Gessel} and $C(n,k)$ is a monic polynomial in $n$, then \begin{align} (2k-1)!! = \sum_{j=0}^k(-1)^{j+k}\binom{2k}{k-j}S(j+k,j). \label{fourth}\end{align}

Compare this result with \[k! = \sum_{j=0}^k(-1)^{j+k}\binom{2k+1}{k-j}S(j+k,j),\] found by setting $n=k+1$ in \eqref{third}. See \cite[equation 14.34]{Gould}. From a survey of identities for the double factorial, we also have \[(2k-1)!! = \sum_{j=0}^k(-2)^{k-j}s(k,j).\] See \cite[identity 5.3]{Callan}. This easily follows by plugging $x = -\frac{1}{2}$ in \eqref{first}. Notice the summation is over non-negative integers whereas identity \eqref{fourth} involves both negative and positive integers. Setting $x=-\frac{1}{2}$ in \eqref{second} gives us \[1=\sum_{k=0}^n(-2)^{n-k}S(n,k)(2k-1)!!,\] an interesting expression but not one that allows us to solve for the odd double factorial like we have in \eqref{fourth}. Identity \eqref{fourth} was initially discovered during the preparation of an expository paper \cite{Hussein} on Wolstenholme's Theorem.




\end{document}